\documentclass[12pt]{article}
\usepackage{amsmath, amscd, amssymb, latexsym,hhline}
\def\<{{\langle}}
\def\>{{\rangle}}
\def\pf{{\it Proof. }}

\def\Alg{{\mbox{\rm Alg}}}

\def\End{{\mbox{\rm End}}}

\def\Tr{{\mbox{\rm Tr}}}

\def\qed{{\quad\mbox{$\blacksquare$} }}

\def\e{{\varepsilon}}
\def\a{{\alpha}}
\def\b{{\beta}}
\def\>d{{\rightharpoondown}}
\def \>u {{\rightharpoonup}}
\def\d<{{\leftharpoondown}}

\def\w{{\omega}}
\def\ba{{\mathbf a}}
\def\bb{{\mathbf b}}
\def\bx{{\mathbf x}}

\def\hbar{{\overline{h}}}

\def\K{{\mathcal K}}

\def\Z{{\bf Z}}

\def\tr{{\mbox{\rm Tr}}}

\numberwithin{equation}{section}

\newtheorem{thm}{Theorem}[section]
\newtheorem{cor}[thm]{Corollary}

\newtheorem{lem}[thm]{Lemma}

\begin{document}
\author{Siu-Hung Ng \\ Mathematics Department, Towson University, Baltimore,
MD 21252}
\title{Hopf Algebras of Dimension $pq$}
\date{}
 \maketitle
\begin{abstract}
Let $H$ be a non-semisimple Hopf algebra with antipode $S$ of dimension $pq$ over
an algebraically closed field of characteristic 0 where $p \le q$
are odd primes. We prove that $\Tr(S^{2p})=p^2d$ where $d \equiv
pq \pmod{4}$. As a consequence,  if $p,q$ are twin primes, then
any Hopf algebra of dimension $pq$ is semisimple.
\end{abstract}

\section{Introduction}
Let $p$ be an odd prime and $k$ an algebraically closed field of
characteristic 0. If $H$ is a semisimple Hopf algebra of dimension
$p^2$ over $k$, then $H$ is isomorphic to
$k[\Z_{p^2}]$ or $k[\Z_p \times \Z_p]$ by \cite{Mas96}. A more general result for
semisimple Hopf algebras of dimension $pq$, where $p$, $q$ are odd
primes, is obtained by \cite{EG99}. In \cite{Ng02},  the author
proved that non-semisimple Hopf algebras of dimension $p^2$ over
$k$ are Taft algebras and hence completed the classification of Hopf
algebras of dimension $p^2$. However, there is no known example
of non-semisimple Hopf algebras of dimension $pq$, with $p<q$. In
fact, it is shown in \cite{AN1} and \cite{BD} that there is no non-semisimple Hopf algebras over $k$ of
dimension 15, 21, 35, 55, 77, 65, 91 or 143.\\

By \cite{Ng02}, if $p\le q$ are odd primes and $H$ is a
non-semisimple Hopf algebra with antipode $S$ of dimension $pq$,
then $S^{4p}=id_H$ and $\Tr(S^{2p})=p^2 d$ for some odd integer
$d$. In this paper, we prove that $d \equiv pq \pmod{4}$. As a
 consequence, we prove that if $p$, $q$ are twin
primes, any Hopf algebra of dimension $pq$ over $k$ is semisimple.
Recently, Etingof and Gelaki
also announce a even more general result \cite{EG03} which covers the cases when
$p < q \le 2p+1$.

\section{Notation and Preliminaries}\label{s1}
Throughout this paper, $p \le q$ are odd primes, $k$ denotes an algebraically closed
 field of characteristic 0, and $H$ denotes a finite-dimensional Hopf algebra over
$k$ with antipode $S$. Its comultiplication and counit are,
respectively, denoted by $\Delta$ and  $\e$. We will use
Sweedler's notation \cite{Sw69}:
$$
\Delta(x) = \sum x_{(1)} \otimes x_{(2)}\,.
$$
A non-zero element $a \in H$ is called group-like if $\Delta(a)=a
\otimes a$.  The set of all group-like elements $G(H)$ of $H$ is a linearly
independent set, and it forms a group under the multiplication of
$H$. For the details of elementary aspects for
finite-dimensional Hopf algebras,  readers are referred to
the references \cite{Sw69} and \cite{Mont93bk}.\\

Let $\lambda \in H^*$ be a non-zero right integral of $H^*$ and  $\Lambda\in H$  a non-zero
 left integral of $H$. There exists $\a \in \Alg(H,k)=G(H^*)$, independent of
the choice of $\Lambda$,  such that $\Lambda a = \a(a) \Lambda$ for  $a \in H$.
Likewise, there is a group-like element $g \in H$, independent of
the choice of $\lambda$,  such that $\b\lambda  = \b(g)\lambda$ for  $\b\in H^*$.
We call $g$ the distinguished group-like element of $H$ and $\a$
the distinguished group-like element of $H^*$. Then we have a
formula for $S^4$ in terms of $\a$ and $g$~\cite{Radf76}:
\begin{equation}\label{eS}
S^4(a) = g(\a \rightharpoonup a \leftharpoonup \a^{-1}) g^{-1}
\quad\mbox{for }a \in H\,,
\end{equation}
where $\rightharpoonup$ and $\leftharpoonup$ denote the natural
actions of the Hopf algebra $H^*$ on $H$ described by
$$
\b\rightharpoonup a = \sum a_{(1)}\b(a_{(2)}) \quad \mbox{and}
\quad a \leftharpoonup \b = \sum \b(a_{(1)})a_{(2)}
$$
for $\b \in H^*$ and $a \in H$.
If $\lambda$ and $\Lambda$ are normalized,
there are formulae for the trace of any linear endomorphism on
$H$.
\begin{thm}{\rm \cite[Theorem 1]{Radf90}}\label{th1.2}
Let $H$ be a finite-dimensional Hopf algebra with antipode $S$ over the field $k$.
Suppose that $\lambda$ is a  right integral of $H^*$, and that $\Lambda$ is
a  left integral of $H$ such that  $\lambda(\Lambda)=1$. Then for any $f \in
\End_k(H)$,
\begin{eqnarray*}
\tr(f)&=&\sum \lambda\left(S(\Lambda_{(2)})f(\Lambda_{(1)}\right))\\
        &=&\sum \lambda\left((S\circ f)(\Lambda_{(2)})\Lambda_{(1)}\right)\\
        &=&\sum \lambda\left((f\circ S)(\Lambda_{(2)})\Lambda_{(1)}\right)\,.
\end{eqnarray*} \hfill\qed
\end{thm}

Following \cite[Section 2]{Ng02}, the index of $H$ is the least
positive integer $n$ such that
$$
S^{4n}=id_H \quad \mbox{and}\quad g^n=1\,.
$$
Suppose that $H$ is a finite-dimensional Hopf algebra of  odd
index $n >1$, and that $\w\in k$ is a primitive $n$th of unity.
Since $g^n=1$ and $\a$ is an algebra map, $\a(g)$ is a $n$th root
of unity. There exists a unique element $x(\w,H) \in \Z_n$ such
that
$$
\a(g)=\w^{x(\w,H)}\,.
$$
Following the notation in \cite{Ng02}, we  let
$$
H^\w_{a, i,j}  =\{ u \in H\,|\, S^2(u) = (-1)^a\w^iu \mbox{ and } ug=\w^j u\}
$$
for any $(a, i, j) \in \Z_2 \times \Z_n \times \Z_n$. Since the $r(g) \in \End_k(H)$, defined by
$r(g)(a) =a g$ for $a \in H$, commutes with $S^2$, we have
\begin{equation}\label{eq:decomposition}
H=\bigoplus_{\ba \in \K_n} H^\w_\ba
\end{equation}
where $\K_n$ denotes the group $\Z_2 \times \Z_n \times \Z_n$.\\

Using the eigenspace decomposition of $H$ in
(\ref{eq:decomposition}), the diagonalization of a left integral
$\Lambda$ of $H$ admits the following form (cf. \cite{Ng02}),
\begin{equation}\label{normalform}
\Delta(\Lambda) = \sum_{\ba \in \K_n} \left(\sum u_\ba \otimes
v_{-\ba+\bx}\right)
\end{equation}
where $\sum u_\ba \otimes v_{-\ba+\bx} \in H_\ba^\w \otimes
H_{-\ba+\bx}^\w$ and $\bx=(0, -x(\w,H), x(\w,H))$ in $\mathcal{K}_n$. \\

In the sequel, we will call the expression in equation
(\ref{normalform}) the {\em normal form} of $\Delta(\Lambda)$
associated with $\w$. We will simply write $u_\ba \otimes
v_{-\ba+\bx}$ for the sum $\sum u_\ba \otimes v_{-\ba+\bx}$ in the
normal form of $\Delta(\Lambda)$.\\

Let $E_\ba^\w$, $\ba \in \K_n$,
be the set of orthogonal projections associated with the decomposition (\ref{eq:decomposition}).
Then
$$
\dim(H_\ba^\w) = \Tr(E_\ba^\w)
$$
 and we have the following lemma.
\begin{lem}\label{l1.1}
Let $H$ be a finite-dimensional Hopf algebra with the antipode $S$
of odd index $n>1$ over $k$, and $\w \in k$  a primitive $n$th
root of unity. Let $x=x(\w,H) \in \Z_n$ and $\bx=(0, -x, x)$.
Then we have
$$
\dim(H^\w_\ba) = \dim(H^\w_{\bx-\ba})
$$
for all $\ba \in \K_n$.
\end{lem}
\pf Let $\Lambda$ be a left integral for $H$ and let $\lambda$ be a right
integral for $H^*$ such that $\lambda(\Lambda)=1$. Using the normal form of
$\Delta(\Lambda)$ associated with $\w$ in (\ref{normalform}) and
Theorem \ref{th1.2}, we have
$$
\Tr(E_\ba^\w)  = \sum_{\bb \in \K_n}
\lambda\left( S(v_{-\bb+\bx})E^\w_{\ba}(u_\bb)\right)=\lambda\left(S(v_{-\ba+\bx})u_\ba\right)\,.
$$
By Theorem \ref{th1.2} again, we also have
$$
\Tr(E_{-\ba+\bx}^\w)  = \sum_{\bb \in \K_n}
\lambda\left( S( E^\w_{-\ba+\bx}(v_{-\bb+\bx}))u_\bb \right)=\lambda\left(S(v_{-\ba+\bx})u_\ba\right)\,.
$$
Therefore, $\Tr(E_\ba^\w)=\Tr(E_{\bx-\ba}^\w)$. Since
 $\dim(H^\w_\ba)  = \Tr(E_\ba^\w)$ for any $\ba \in \K_n$, the result follows. \qed \\

\begin{thm}{\rm \cite{Ng02}} \label{t:2.2}
Let $H$ be a  Hopf algebra of dimension
$pq$ over $k$ with antipode $S$, where $p \le q$ are odd primes.
Then the index of $H$ and the order of $S^4$ are equal to $p$,
and $\Tr(S^{2p})= p^2 d$ for some odd integer $d$.\\
\mbox{}\hfill\qed\\
\end{thm}

\begin{lem}\label{l:2.4}
Suppose that $H$ is a non-semisimple Hopf algebra of dimension
$pq$ over $k$ where $p \le q$ are odd primes,  and that $\w \in k$
is a primitive $p$th root of unity. Let $g$ and $\a$ be the
distinguished group-like elements of $H$ and $H^*$ respectively.
If $g$ is non-trivial, then the integer $d$ in Theorem \ref{t:2.2}
is given by
$$
\dim(H^\w_{0, i, j}) - \dim(H^\w_{1, i, j}) =d
$$
for all $i, j \in \Z_p$. Moreover, if both $g$ and $\a$ are not
trivial, then
$$
\dim(H^\w_{a, i, j'} )= \dim(H^\w_{a, i, j})
$$
for any $a \in \Z_2$ and $i, j, j' \in \Z_p$.
\end{lem}
\pf
If $\a$ is trivial and $g \ne 1$, then by \cite[Lemma 4.3]{Ng02},
$$
\dim(H^\w_{0, i, j}) - \dim(H^\w_{1, i, j})=d\,.
$$
If both $g$ and $\a$ are non-trivial, then by the proof of
\cite[Proposition 5.3]{Ng02},  $H$ is isomorphic to the biproduct
\begin{equation} \label{e:byproduct}
R \times B
\end{equation}
as Hopf algebras (cf. \cite{Radf85}) where $B=k[g]$  and $R$ is a
right $B$-comodule subalgebra of $H$. It is shown in \cite[section
4]{AnSc98} that $R$ is invariant under $S^2$. Moreover, in the
identification $H \cong R \otimes B$ given by multiplication, one
has
\begin{equation}\label{S_square}
S^2 = T \otimes id_B
\end{equation}
where $T$ is the restriction of $S^2$ on $R$. Let
$$
R_{a,i} = \{ x \in R \,|\, S^2(x) = (-1)^a \w^i x\}
$$
for any $(a, i) \in \Z_2 \times \Z_p$. It follows from the proof
of \cite[Proposition 5.3]{Ng02} that
$$
\dim(R_{0, i}) - \dim(R_{1, i}) = d\,.
$$
By (\ref{e:byproduct}),
$$
H^\w_{a, i, j} = R_{a,i} \otimes e_j
$$
for all $(a, i, j) \in \mathcal{K}_p$ where $e_j$ is the central idempotent of $B$ such that
$e_j g= \w^j e_j$. Thus,
$$
\dim(H^\w_{a, i, j}) = R_{a,i}
$$
for all $(a, i,j) \in  \mathcal{K}_p$ and hence
$$
\dim(H^\w_{0, i, j}) - \dim(H^\w_{1, i, j}) = d\,.
$$
\mbox{}\hfill\qed

\section{Proofs of Main Results}
\begin{lem} \label{l1}
Let $H$ be a finite-dimensional Hopf algebra with antipode $S$ of
odd index $n>1$ over $k$, and $\w \in k$  a primitive $n$th root
of unity. Let $\ell \in \Z_n$ such that $2\ell=x(\w,H)$. Then
$$
\dim(H^\w_{1,-\ell, \ell})
$$
is even.
\end{lem}
\pf Let $V$ be space of all $f \in H^*$ such that $f(u)=0$ for $u \in
H^\w_{a,i,j}$ whenever $(a,i,j) \ne (1, -\ell,\ell)$.
Obviously, $V$ is isomorphic to $(H^\w_{1,-\ell, \ell})^*$ and so
$\dim(V) = \dim(H^\w_{1, -\ell, \ell})$.
Let $\Lambda$ be a non-zero left integral of $H$ and
$$
\Delta(\Lambda) = \sum_{\ba \in \K_n} u_\ba \otimes
v_{-\ba+\bx}
$$
the normal form of $\Delta(\Lambda)$ associated with $\w$ where
$\bx=(0, -2\ell, 2\ell)$. Then
$$
(f,h) = (f \otimes h)\Delta(\Lambda)
$$
defines a non-degenerate bilinear form on $H^*$. Let $f \in V$
such that $(f,h)=0$ for all $h \in V$. For any  $h' \in H^*$,
there exists $h \in V$ such that $h'(u)=h(u)$ for all $u \in
H^\w_{1,-\ell, \ell}$. Thus
$$
(f, h')=\sum_{\ba \in \K_n} f(u_\ba)h'(v_{-\ba+\bx}) = f(u_{1, -\ell, \ell})h'(v_{1, -\ell, \ell})
= (f, h)=0\,.
$$
By the non-degeneracy of $(\cdot, \cdot)$, $f=0$. Therefore, $(\cdot, \cdot)$ induces a non-degenerate
bilinear form on $V$. Using
\cite[Theorem 3(d)]{Radf94}, we have
$$
\Delta^{op}(\Lambda) = \sum_{(a,i,j) \in \K_n} (-1)^a
\w^{-i-j}\left(\sum u_{a,i,j} \otimes
v_{a,-2\ell-i,2\ell-j}\right)\,.
$$
Therefore, for any $f,h \in V$,
$$
(h,f) = (f \otimes h)\Delta^{op}(\Lambda)=-f(u_{1,-\ell, \ell})h(v_{1,-\ell, \ell})=-(f,h)\,.
$$
Hence, $V$ admits a non-degenerate alternating form and so
$\dim(V)$ is even. \qed\\

If $H$ is a finite-dimensional Hopf algebra of index $n>1$, we
define
\begin{eqnarray*}
   H_-&:=&\{u \in H\,|\, S^{2n}(u)=-u\}\,,\\
   H_+&:=&\{u \in H\,|\, S^{2n}(u)=u\}\,.\\
\end{eqnarray*}

\begin{cor} \label{l:3.2}
Suppose $H$ is a finite-dimensional Hopf algebra with antipode $S$
of odd index $n>1$. Then, the subspace $H_-$ is of even dimension.
\end{cor}
\pf Let $\w \in k$ be a $n$th of unity and $\ell \in \Z_n$  such
that $2\ell = x(\w,H)$. We then have
$$
  H_-  =  \bigoplus_{i,j \in \Z_n} H^\w_{1,i,j}
       =  H^\w_{1, -\ell, \ell}
\oplus \left(\bigoplus_{ {\rm some}\,\,i,j \in \Z_n \atop
                                      (i,j) \ne (-\ell, \ell)}
H^\w_{1,i,j} \oplus H^\w_{1,-2\ell-i,2\ell-j}\right)\,.
$$
It follows from Corollary \ref{l1.1} and Lemma
\ref{l1}, $\dim(H_-)$ is even. \qed\\

\begin{thm}\label{t:d}
Let $H$ be a non-semisimple Hopf algebra with antipode $S$ of
dimension $pq$ where $p\le q$ are odd primes. Then
$$
\Tr(S^{2p})=p^2d \quad\mbox{and} \quad d \equiv pq \pmod{4}\,.
$$
\end{thm}
\pf By Theorem \ref{t:2.2}, $H$ is of index $p$ and $\Tr(S^{2p}) =
p^2d$ for some odd integer $d$. Since
$$
\dim(H_+) + \dim(H_-) = pq
$$
 and
 $$
 \Tr(S^{2p}) = \dim(H_+) -\dim(H_-) = p^2d\,,
 $$
we have
$$
\dim(H_-) = p (q-pd)/2\,.
$$
By Corollary \ref{l:3.2}, $p(q-pd) \equiv 0 \pmod{4}$ or $d \equiv
pq \pmod{4}$\,. \qed\\
\begin{thm}\label{t:final}
For any pair of twin primes $p<q$, if $H$ is a Hopf algebra of dimension $pq$, then
$H$ is semisimple.
\end{thm}
\pf Suppose there is a non-semisimple Hopf algebra $H$ of dimension
$pq$. By \cite{LaRa88}, $H^*$ is also non-semisimple. Since $\dim(H)$ is odd, by
\cite[Theorem 2.1]{LaRad95}, $H$ and $H^*$ can not be both unimodular.
By duality, we may simply assume that $H^*$ is not unimodular. It follows from
Theorem \ref{t:2.2} that $|G(H)|=p$ and so
$$
\dim(C) \ge p
$$
where $C$ is the coradical of $H$. If $\dim(C) = p$, then $H$ is pointed and hence, by \cite[Corollary 4]{Stefan97},
$H$ is semisimple. Therefore, $\dim(C) > p$ and so we have
$$
\Tr(S^{2p}|_{H/C}) \ge -(pq-\dim(C)) > -pq + p = -p^2 -p\,.
$$
It follows from \cite[Lemma 3.2]{LaRa88} that
$$
\Tr(S^{2p}|_C) \ge p\,.
$$
Thus, we have
\begin{equation}\label{eq:trace}
\Tr(S^{2p}) = \Tr(S^{2p}|_C) + \Tr(S^{2p}|_{H/C})> -p^2\,.
\end{equation}
Since $pq \equiv -1 \pmod{4}$, by Theorem \ref{t:d},
$$
\Tr(S^{2p})=-p^2
$$
but this contradicts (\ref{eq:trace}). \qed

\begin{center}
{\bf Acknowledgement}
\end{center}
The author would like thank P. Etingof for his useful suggestion for Theorem
\ref{t:final} and informing me his recent work \cite{EG03} with S. Gelaki.
\bibliographystyle{amsalpha}

\begin{thebibliography}{Mon93}

\bibitem[AN01]{AN1}
Nicol{\'a}s Andruskiewitsch and Sonia Natale, \emph{Counting arguments for
  {H}opf algebras of low dimension}, Tsukuba J. Math. \textbf{25} (2001),
  no.~1, 187--201. \MR{2002d:16046}

\bibitem[AS98]{AnSc98}
Nicol{\'a}s Andruskiewitsch and Hans-J{\"u}rgen Schneider, \emph{Hopf algebras
  of order $p\sp 2$ and braided {H}opf algebras of order $p$}, J. Algebra
  \textbf{199} (1998), no.~2, 430--454. \MR{99c:16033}
\bibitem[BD02]{BD}
Margaret Beattie and Sorin Dascalescu, \emph{Hopf Algebras of Dimension 14},
Preprint {\bf arXiv:math.QA/0205243}
\bibitem[EG98]{EG99}
Pavel Etingof and Shlomo Gelaki, \emph{Semisimple {H}opf algebras of dimension
  $pq$ are trivial}, J. Algebra \textbf{210} (1998), no.~2, 664--669.
  \MR{99k:16079}
\bibitem[EG03]{EG03}
Pavel Etingof and Shlomo Gelaki, \emph{On Hopf Algebras of Dimension $pq$},
Preprint {\bf arXiv:math.QA/0303359}

\bibitem[LR88]{LaRa88}
Richard~G. Larson and David~E. Radford, \emph{Finite-dimensional cosemisimple
  {H}opf algebras in characteristic $0$ are semisimple}, J. Algebra
  \textbf{117} (1988), no.~2, 267--289. \MR{89k:16016}

\bibitem[LR95]{LaRad95}
\bysame, \emph{Semisimple {H}opf algebras}, J. Algebra \textbf{171} (1995),
  no.~1, 5--35. \MR{96a:16040}

\bibitem[Mas96]{Mas96}
Akira Masuoka, \emph{The $p\sp n$ theorem for semisimple {H}opf algebras},
  Proc. Amer. Math. Soc. \textbf{124} (1996), no.~3, 735--737. \MR{96f:16046}

\bibitem[Mon93]{Mont93bk}
Susan Montgomery, \emph{Hopf algebras and their actions on rings}, CBMS
  Regional Conference Series in Mathematics, vol.~82, Published for the
  Conference Board of the Mathematical Sciences, Washington, DC, 1993.

\bibitem[Ng02]{Ng02}
Siu-Hung Ng, \emph{Non-semisimple Hopf algebras of dimension $p\sp 2$}, J.
  Algebra \textbf{255} (2002), no.~1, 182--197.

\bibitem[Rad76]{Radf76}
David~E. Radford, \emph{The order of the antipode of a finite dimensional
  {H}opf algebra is finite}, Amer. J. Math. \textbf{98} (1976), no.~2,
  333--355. \MR{53 \#10852}

\bibitem[Rad85]{Radf85}
\bysame, \emph{The structure of {H}opf algebras with a projection}, J. Algebra
  \textbf{92} (1985), no.~2, 322--347. \MR{86k:16004}

\bibitem[Rad90]{Radf90}
\bysame, \emph{The group of automorphisms of a semisimple {H}opf algebra over a
  field of characteristic $0$ is finite}, Amer. J. Math. \textbf{112} (1990),
  no.~2, 331--357. \MR{91b:16048}

\bibitem[Rad94]{Radf94}
\bysame, \emph{The trace function and {H}opf algebras}, J. Algebra \textbf{163}
  (1994), no.~3, 583--622. \MR{95e:16039}
\bibitem[{\c{S}}te97]{Stefan97}
D.~{\c{S}}tefan, \emph{Hopf subalgebras of pointed {H}opf algebras and
  applications}, Proc. Amer. Math. Soc. \textbf{125} (1997), no.~11,
  3191--3193. \MR{97m:16076}

\bibitem[Swe69]{Sw69}
Moss~E. Sweedler, \emph{Hopf algebras}, W. A. Benjamin, Inc., New York, 1969,
  Mathematics Lecture Note Series.

\end{thebibliography}
\providecommand{\bysame}{\leavevmode\hbox to3em{\hrulefill}\thinspace}
\providecommand{\MR}{\relax\ifhmode\unskip\space\fi MR }
\providecommand{\MRhref}[2]{%
  \href{http://www.ams.org/mathscinet-getitem?mr=#1}{#2}
}
\providecommand{\href}[2]{#2}

\end{document}